\documentclass[11pt]{amsart}
\usepackage[english]{babel}
\usepackage{amssymb}
\usepackage{fourier}
\usepackage{mathrsfs}
\usepackage{enumerate}
\usepackage{color}
\usepackage{ifpdf}
\usepackage{tikz}
\usepackage{tikz-cd}
\usetikzlibrary{decorations.pathmorphing,arrows}
\usepackage[initials]{amsrefs}

\usepackage{caption}
\usepackage{subcaption}
\usepackage{comment}
\usetikzlibrary{intersections}

\newcommand{\bbN}{{\mathbb N}}
\newcommand{\bbQ}{{\mathbb Q}}
\newcommand{\bbR}{{\mathbb R}}

\newcommand{\bbZ}{{\mathbb Z}}
\newcommand{\bbC}{{\mathbb C}}

\newcommand{\bfH}{{\mathbf H}}

\newcommand{\bfZ}{{\mathbf Z}}

\newcommand{\PGL}{\operatorname{PGL}}
\newcommand{\SL}{\operatorname{SL}}
\newcommand{\PSL}{\operatorname{PSL}}

\newcommand{\setdef}[2]{ \left\{ {#1}\ \left/\ {#2} \right.\right\} }

\newtheorem{theorem}{Theorem}[section]
\newtheorem{thm}[theorem]{Theorem}
\newtheorem{mthm}{Theorem}

\newtheorem{lemma}[theorem]{Lemma}

\newtheorem{corollary}[theorem]{Corollary}

\newtheorem{conjecture}[theorem]{Conjecture}

\newtheorem{claim}[theorem]{Claim}

\theoremstyle{definition}
\newtheorem{definition}[theorem]{Definition}

\newtheorem{remark}[theorem]{Remark}

\numberwithin{equation}{section}

\title[B-C property and arithmetic Kleinian groups]{Bounded clustering property characterizes arithmetic nonuniform Kleinian groups}
\author{Yanlong Hao}
\address{University of Illinois at Chicago}
\email{yhao23@uic.edu}

\begin{document}
\date{\today}

\begin{abstract}
In this short note, we show that the B-C property characterizes arithmetic lattices along all nonuniform lattices in $\PSL(2,\bbC)$.
\end{abstract}

\maketitle

\section{Introduction and Statement of the main results} 
For all $m$, $n\in \bbZ$, denote $$S(m,n)=\setdef{z\in \mathbb{C}}{m\leq{\text{Re}}(z)\leq m+1, n\leq {\text{Im}}(z)\leq n+1}.$$
We say that a set $A$ of complex numbers satisfies the \emph{bounded clustering} or \emph{B-C property}
iff there exists a constant $K_A$ such that $A\cap S(m,n)$ has less than $K_A$ elements for all $m$, $n\in \bbZ$. Further set
$${\text{Gap}}(A):=\inf\setdef{|a-b|}{a,b\in A, a\neq b}$$

Let $\Gamma$ be a Kleinian group, i.e., a discrete subgroup of $\PSL(2,\mathbb{C})$. The \emph{trace set} $\text{Tr}(\Gamma)$ of $\Gamma$ is defined (up to a sign) as the set of traces of elements of $\Gamma$.

In \cite{luo1994number}, Luo and Sarnak showed the trace set of an arithmetic Fuchsian group satisfies the B-C property. The proof indeed also works for Kleinian groups. Furthermore, Sarnak conjectured that the converse is also true.
\begin{conjecture}[Sarnak \cite{MR1321639}]\label{BC}
   Let $\Gamma$ be a cofinite Fuchsian group. 
   \begin{enumerate}
       \item 
   If $\mathrm{Tr}(\Gamma)$ satisfies the B-C property, then $\Gamma$ is arithmetic. 
   \item
   If $\mathrm{Gap(Tr}(\Gamma))>0$, then $\Gamma$ is derived from a quaternion algebra.
   \end{enumerate}
\end{conjecture}

In \cite{MR1394753}, Schmutz makes an even stronger conjecture using the linear growth of a set instead of the B-C property. A subset of reals is said to have linear growth if and only if there exist positive real constants $C$ and $D$ such that for all $n$,
$$\#\setdef{a\in A}{|a|\leq n}\leq Cn+D.$$
\begin{conjecture}[Schmutz \cite{MR1394753}]\label{LG}
    Let $\Gamma$ be a cofinite Fuchsian group. If
${\mathrm{Tr}}(\Gamma)$ has linear growth then $\Gamma$ is arithmetic.
\end{conjecture}
 Schmutz proposed a proof of Conjecture~\ref{LG} for nonuniform lattices. In this paper, Schmutz essentially proved part (2) of Sarnak's Conjecture~\ref{BC} under (1). Unfortunately, the proof of Conjecture~\ref{LG} contains a gap. Later, Geninska and Leuzinger fixed this gap in \cite{geninska2008geometric} and confirmed part (1) of Sarnak's Conjecture for nonuniform Fuchsian lattices.  Note that Conjecture~\ref{LG} is still open even for nonuniform lattice. And the Conjecture~\ref{BC} remains entirely open for cocompact Fuchsian groups.

In this paper, we generalize the work of Geninska and Leuzinger to Kleinian groups. And give a pure algebraic proof of Conjecture~\ref{BC} in the nonuniform case.
\begin{mthm}\label{main}
Let $\Gamma$ be a non-uniform lattice of  $\PSL(2,\mathbb{R})$ or $\PSL(2,\mathbb{C})$. 
   \begin{enumerate}
       \item 
   If $\mathrm{Tr}(\Gamma)$ satisfies the B-C property, then $\Gamma$ is arithmetic. 
   \item
   If $\mathrm{Gap(Tr}(\Gamma))>0$, then $\Gamma$ is derived from a quaternion algebra.
   \end{enumerate}
\end{mthm}

The strategy of proving Theorem~\ref{main} is similar to Schmutz's proposal. Given an element $\gamma$ in $\Gamma$, Schmutz constructs a $Y$-piece $S$, a surface of signature (0,3) related to $\gamma$. By considering traces of different families of elements in $\pi_1(S)$, there are  restrictions on the trace of $\gamma$. The whole process could be translated into pure algebraic language. However, the algebraic version representation of those elements involved in Schmutz's work is quite complicated. Also, it is unclear whether the same construction outcomes for Klein groups. 

Our approach relies heavily on unipotent elements. Two different cusp subgroups provide a few families of elements. It turns out that the traces of these elements are enough for us to have the same conclusion as in the work of Schnutz,  Geninska, and Leuzinger.

The note is organized as follows: Section 2 recalls some preliminaries on (arithmetic) Fuchsian and Kleinian groups. Then in section 3, we prove Theorem~\ref{main} for Fuchsian groups. Most of the proof is inspired by \cite{geninska2008geometric}. In section 4, a similar strategy extends to Kleinian groups. Finally, in section 5, there is a corollary for personal interest. 

\medskip

I want to acknowledge and thank Alexander Furman for his many suggestions and Tian Wang for discussions of number theory.  

\section{Definitions, notations, and some preliminaries}
A general reference for this section is the book \cite{maclachlan2003arithmetic}.
\subsection{Fuchsian and Kleinian groups}

We denote by $\SL(2,\mathbb{R})$
the group of real $2\times2$ matrices with determinant 1 and by $\PSL(2,\mathbb{R})$ the
quotient group $\SL(2,\mathbb{R}))/\{\pm I_2\}$ where $I_2$ is the $2\times 2$ identity matrix. Similarly, we have $\PSL(2,\mathbb{C})$.

A discrete subgroup of $\PSL(2,\mathbb{R})$ and $\PSL(2,\mathbb{C})$ is called \emph{Fuchsian} and \emph{Kleinian group}, respectively. 

Let $p_F:\SL(2,F)\rightarrow\PSL(2,F)$ be the projection where $F=\mathbb{R}$ or $\mathbb{C}$. Give the Fuchsian or Klein group $\Gamma$, denote $\bar{\Gamma}=p_F^{-1}(\Gamma)$. We call
$$\mathrm{Tr}(\Gamma):=\setdef{\mathrm{tr}T}{T\in \bar{\Gamma}}$$
the \emph{trace set} of $\Gamma$. And similarly,
\emph{reduce trace set} of $\Gamma$ is 
$$\mathrm{Tr_R}(\Gamma):=\setdef{\mathrm{tr}T}{T\in \bar{\Gamma}, T\neq \pm I_2}$$

A lattice of a locally compact, second countable topological group $G$ is a discrete subgroup $\Gamma$ such that $G/\Gamma$ has finite Haar measure. A lattice is called uniform if $G/\Gamma$ is compact and nonuniform otherwise. 

A Fuchsian or Kleinian lattice $\Gamma$ is nonuniform iff $\Gamma$ contains parabolic elements. 

\medskip 

Arithmetic Kleinian groups are obtained in the following way: Let $k$ be an algebraic number field
with exactly one non-real archimedean place so that the $\bbQ$-isomorphisms of $k$ into
$\mathbb{C}$ are $\phi_1$, $\phi_2$, $\ldots$, $\phi_n$ where we take $\phi_2=\mathrm{Id}$, $\phi_2$ the complex conjugation and $\phi_i(k)\subset \mathbb{R}$ for $i=3,4,...,n$. Let $A$ be a quaternion algebra over $k$, which is ramified at all real places, and thus there is an isomorphism
$$\rho: A\otimes_\bbQ\mathbb{R}=\SL(2,\mathbb{C})\oplus \bfH\oplus \bfH\oplus\cdots\oplus\bfH,$$
where $\bfH$ denotes Hamilton's quaternions. Denote $P$ the projection to the first factor.

Let $\mathcal{O}$ be an order in $A$, and $\mathcal{O}^1$  denote the group of elements of reduced norm 1. Then $P\rho(\mathcal{O}^1)$ is a lattice of $\PSL(2,\mathbb{C})$. The class of
\emph{arithmetic Kleinian groups} is all Kleinian lattices commensurable with such groups $P\rho(\mathcal{O}^1)$. In addition, we say that a Kleinian group is \emph{derived from a
quaternion algebra} if it is a subgroup of finite index in some $P\rho(\mathcal{O}^1)$.

Arithmetic Fuchsian groups are similarly defined: In this case, the field
$k$ is to be totally real, and the quaternion algebra $A$ is ramified at all real places except one, which can be taken to be the identity. There is thus an isomorphism
$$\rho: A\otimes_\bbQ\mathbb{R}=\SL(2,\mathbb{R})\oplus \bfH\oplus \bfH\oplus\cdots\oplus\bfH.$$
One obtains arithmetic Fuchsian groups, Fuchsian groups derived from quaternion algebras precisely as above.
\subsection{Characterization of arithmetic Fuchsian and Kleinian groups}
Takeuchi characterizes arithmetic Fuchsian groups in the class of all Fuchsian lattices in \cite{takeuchi1975characterization}. Maclachlan and Reid \cite{maclachlan1987commensurability} generated this work for Kleinian groups. 

Let $\Gamma$ be a Fuchsian or Kleinian group, and let $\Gamma^{(2)}$ denote the subgroup generated by the
squares of elements of $\Gamma$. Note that if $\Gamma$ is finitely generated then $\Gamma^{(2)}$ is of finite index
in $\Gamma$.
\begin{thm}[\cite{takeuchi1975characterization}, \cite{borel1981commensurability}]\label{arithmetic}
    If $\Gamma$ is an arithmetic Fuchsian or Kleinian group, then $\Gamma^{(2)}$ is derived from a quaternion algebra.
\end{thm}
\begin{thm}[\cite{takeuchi1975characterization}]\label{THM: Fuchsian}
Let $\Gamma$ be a cofinite Fuchsian group. Then $\Gamma$ is derived from a quaternion algebra over a totally real algebraic number field if and
only if $\Gamma$ satisfies the following two conditions:
\begin{enumerate}
    \item 
$K:=\bbQ(\mathrm{Tr}(\Gamma))$ is an algebraic number field of finite degree and $\mathrm{Tr}(\Gamma)$
is contained in the ring of integers $\mathcal{O}_K$ of $K$.
\item
For any embedding $\phi$ of $K$ into $\mathbb{C}$, which is not the identity, $\phi(\mathrm{Tr}(\Gamma))$
is bounded in $\mathbb{C}$.
\end{enumerate}
\end{thm}
\begin{thm}[\cite{maclachlan1987commensurability}]\label{THM: Kleinian}
Let $\Gamma$ be a cofinite Kleinian group. Then $\Gamma$ is derived from a quaternion algebra if and
only if $\Gamma$ satisfies the following two conditions:
\begin{enumerate}
    \item 
$K:=\bbQ(\mathrm{Tr}(\Gamma))$ is an algebraic number field of finite degree and $\mathrm{Tr}(\Gamma)$
is contained in the ring of integers $\mathcal{O}_K$ of $K$ and $K\nsubseteq \mathbb{R}$.
\item
For any embedding $\phi$ of $K$ into $\mathbb{C}$, which is not the identity, complex conjugation, $\phi(\mathrm{Tr}(\Gamma))$
is bounded in $\mathbb{C}$.
\end{enumerate}
\end{thm}
\section{Fuchsian groups}
We prove Theorem~\ref{main} for Fuchsian groups in this section. 

 Let $\Gamma$ be a non-uniform Fuchsian group and $x\in \partial\Gamma$ be a cusp point with a cusp subgroup $\Gamma_x$. Taking $g\in \Gamma$ with $gx\neq x$, then $gx$ is a cusp point of $\Gamma$, and $\Gamma_{gx}=g\Gamma_x g^{-1}$. Up to conjugation in $\PSL(2,\mathbb{R})$, we may assume $x=\infty$, $gx=0$, and $$\Gamma_\infty=\setdef{\begin{pmatrix}
    1&k\\
    0&1
\end{pmatrix}}{k\in \bfZ}.$$

Since $g\infty=0$, $g$ is in the form $\begin{pmatrix}
    0&\frac{1}{\beta}\\
    -\beta&*
\end{pmatrix}$ for some $\beta\in\bbR$. It follows from the fact $\Gamma_0=g\Gamma_\infty g^{-1}$ that
$$\Gamma_0=\setdef{\begin{pmatrix}
    1&0\\
    k\beta^2&1
\end{pmatrix}}{k\in \bfZ}.$$

From now on, we assume $\Gamma$ is a Fuchsian lattice containing $\Gamma_0$ and $\Gamma_\infty$ as above in this section. We will show that with this embedding, $\Gamma$ commensurate to $\PSL(2,\bbZ)$. To prove it, we first show that $\Gamma^{(2)}<\PSL(2,\bbQ)$ by linear growth of the trace set. Then using the fact trace set has B-C property, the trace set of $\Gamma^{(2)}$ is a subset of $\bbZ$.

\subsection{First step: $\Gamma^{(2)}$ is rational.}
The main result in this subsection is Lemma~\ref{A^2}. 

 Let $x_1,x_2,\cdots,x_n\in\bbR$. Denote $\bbQ\langle x_1,x_2,\cdots,x_n\rangle$ the $\bbQ$-vector space generated by $x_1$, $x_2,\cdots, x_n$. 
\begin{lemma}\label{A^2}
If $\mathrm{Tr}(\Gamma)$ has linear growth, then $A^2\in \PSL(2,\bbQ)$ for all $A\in \Gamma$.
\end{lemma}
\begin{proof}
Let $A=\begin{pmatrix}
    a&b\\c&d
\end{pmatrix}\in \Gamma$.
Without loss of generality, we assume $c\neq 0$. Indeed, $c=0$ implies that $A\infty=\infty$. Then $A\in \Gamma_\infty\subset \PSL(2,\bbQ)$.

Since
$$\begin{pmatrix}
  1&0\\
  k\beta^2&1
\end{pmatrix}\begin{pmatrix}
    a&b\\
    c&d
\end{pmatrix}\begin{pmatrix}
    1&l\\
    0&1
\end{pmatrix}=\begin{pmatrix}
    a&al+b\\
    k\beta^2a+c&kl\beta^2a+lc+k\beta^2b+d
\end{pmatrix},$$
$\mathrm{Tr}(\Gamma)$ contains all elements of the from $a+d+kl\beta^2a+k\beta^2b+lc$, $k,l\in \bbZ$. Hence the set
\[
    \Omega_{a,b,c}:=\setdef{kl\beta^2a+k\beta^2b+lc}{k,l\in\bbZ}
\]
has linear growth. For convenience, denote $\Theta(k,l)=kl\beta^2a+k\beta^2b+lc$.
\begin{claim}\label{rank 1}
    $\bbQ\langle\beta^2a,\beta^2b,c,\beta^2d\rangle=\bbQ\langle c\rangle.$
\end{claim}
The proof of Claim~\ref{rank 1} has two steps.

\textbf{Step 1.} $\beta^2 a\in \bbQ\langle c,\beta^2b\rangle.$

Assume $\beta^2 a\notin \bbQ\langle c,\beta^2b\rangle.$ On one hand, under this assumption, $\Theta(k,l)=\Theta(k',l')$ iff $kl=k'l'$ and $k\beta^2 b+lc=k'\beta^2 b+l'c.$ Assume that $kl\neq 0$, then $l'=\frac{kl}{k'}$, we have
$$(k-k')\beta^2 b=(\frac{l(k-k')}{k'})c.$$

If $k=k'$, then $l=l'$. If $k\neq k'$, then $\frac{l}{k'}=\frac{\beta^2 b}{c}$. Therefore $\Theta$ is at most 2 to 1 when $kl\neq 0$.

On the other hand, the set 
$$D_N:=\setdef{(k,l)}{1\leq kl\leq N,k,l\in \bbN^+}$$
has at least $N\ln N-N$ many elements. And all elements in $\Theta(D_N)$ have absolute value less that $N(|\beta^2 a|+|\beta^2 b|+|c|)$. It is contrary to the linear growth property.
Therefore $\beta^2 a\in \bbQ\langle c,\beta^2b\rangle.$

\textbf{Step 2.} $\beta^2 b\in \bbQ\langle c\rangle.$ 

Since $\beta^2 a\in \bbQ\langle c,\beta^2b\rangle,$ there exist $s,t\in\bbQ$ with $\beta^2 a=s\beta^2 b+t c$. Now 
\[\Omega_{a,b,c}=\setdef{(skl+k)\beta^2b+(tkl+l)c}{k,l\in\bbZ}\]
has linear growth. 

 Define $$\Phi(k,l)=(skl+k,tkl+l).$$ And notice that $$\Theta(k,l)=(skl+k)\beta^2b+(tkl+l)c.$$

Case 1: $s=t=0$. $\Theta$ map the set $\setdef{(k,l)}{1\leq l,k\leq N}$ to numbers with norm no more that $N(|\beta^2 b|+|c|)$. 

Case 2: $s=0$, $t\neq 0$.
 The image of $\Phi$ determines $k$. It follows that $\phi$ is injective when $k\neq -\frac{1}{t}$. The set $D_N\setminus \{k=-\frac{1}{t}\}$ has more than $N\ln N-2N$ elements. On the other hand, $|\Theta(u)|\leq N[(|s|+1)|\beta^2 b|+(|t|+1)|c|]$ for all $u\in D_N\setminus \{k=-\frac{1}{t}\}$. 

Case 3: $s\neq 0$ and $t=0$. Similar to Case 2.

Case 4: $st\neq 0$. On one hand, $|\Theta(u)|\leq N[(|s|+1)|\beta^2 b|+(|t|+1)|c|]$ for all $u\in D_N$. On the other hand, $\Phi$ is at most 2 to 1. 

In all cases, Dirichlet‘s principle gives  $(k,l)\neq(k',l')$ such that $\Theta(k,l)=\Theta(k',l')$ and $\Phi(k,l)\neq\Phi(k',l')$. We have a nontrivial homogeneous linear equation of $\beta^2b$ and $c$. By assumption, $c\neq 0$. Hence $\beta^2 b\in \bbQ\langle c\rangle.$

Similarly, $\beta^2 d\in \bbQ\langle c\rangle$ by
\[
\begin{pmatrix}
    1&l\\
    0&1
\end{pmatrix}\begin{pmatrix}
    a&b\\
    c&d
\end{pmatrix}\begin{pmatrix}
    1&0\\
    k\beta^2&1
\end{pmatrix}=\begin{pmatrix}
    a+kl\beta^2d+lc+k\beta^2b&ld+b\\
    k\beta^2d+c&d
\end{pmatrix}.
\]

The claim is proved.

Now we prove the lemma.

Considering the element $\begin{pmatrix}
    1&1\\
    0&1
\end{pmatrix}\begin{pmatrix}
    1&0\\
    \beta^2&1
\end{pmatrix}=\begin{pmatrix}
    1+\beta^2&1\\
    \beta^2&1
\end{pmatrix}\in\Gamma$.
Claim~\ref{rank 1} gives $\beta^2+\beta^4\in \bbQ\langle \beta^2\rangle$. We conclude that $\beta^2\in \bbQ$. Then 
\[\bbQ\langle a,b,c,d\rangle=\bbQ\langle c\rangle.\]

Now $A=cA'$ with $A'\in \PGL(2,\bbQ)$. Taking determinant, $c^2\in \bbQ$. Finally,
\[A^2=c^2A'^2\in \PSL(2,\bbQ).\] 
\end{proof}

\subsection{Step two: $\mathrm{Tr}(\bar{\Gamma})$ is in $\bbZ$.}

By Theorem~\ref{arithmetic} and Theorem~\ref{THM: Fuchsian}, it is enough to work with $\Gamma^{(2)}$ to show that $\Gamma$ is arithmetic. But to show that $\Gamma$ is derived from a quaternion algebra, we will need a slightly bigger subgroup $\bar{\Gamma}=\Gamma\cap \PSL(2,\bbQ)$.  
 $\bar{\Gamma}$ is a finite index subgroup of $\Gamma$ since $\Gamma$ is finitely generated and $\Gamma^{(2)}<\bar{\Gamma}$.  
 \begin{lemma}\label{Z}
 If $\mathrm{Tr}(\bar{\Gamma})$ has the B-C property, then $\bar{\Gamma}$ is derived from a quaternion algebra.
 \end{lemma}
 \begin{proof}
 Choose $A=\begin{pmatrix}
    a&b\\
    c&d
\end{pmatrix}\in \bar{\Gamma}.$ Then
$$
\begin{pmatrix}
    1&1\\
    0&1
\end{pmatrix}\begin{pmatrix}
    a&b\\
    c&d
\end{pmatrix}\begin{pmatrix}
    1&-1\\
    0&1
\end{pmatrix}\begin{pmatrix}
    d&-b\\
    -c&a
\end{pmatrix}=\begin{pmatrix}
    1+ac+c^2&1-ac-a^2\\
    c^2&1-ac
\end{pmatrix}\in \bar{\Gamma}.$$
By iterating, we get a sequence of elements $A_n\in \bar{\Gamma}$ such that
$A_n=\begin{pmatrix}
   *&*\\
   c^{2^n}&*
\end{pmatrix}$ and $\mathrm{tr}(A_n)=2+c^{2^n}$.

It is clear that $\mathrm{tr}(A_n\begin{pmatrix}
    1&k\\
    0&1
\end{pmatrix})=2+(k+1)c^{2^n}$.
It implies that the set
$$\Delta_c:=\setdef{kc^{2^n}}{k\in \bbZ, n\in\bbN^+}$$
satisfies the B-C property.

In proof of Proposition 4.9 in \cite{geninska2008geometric}, it was shown that when $c\notin \bbZ$, $\Delta_c$ has no B-C property. We conclude $c\in \bbZ$. The same result of a slightly complicated case of Kleinian groups will be proven in the next section in full detail.

Considering the following two equations:
$$\begin{pmatrix}
    1&0\\
    \beta^2&1
\end{pmatrix}\begin{pmatrix}
    a&b\\
    c&d
\end{pmatrix}=\begin{pmatrix}
    a&b\\
    \beta^2 a+c&\beta^2 b+d
\end{pmatrix},$$
$$\begin{pmatrix}
    a&b\\
    c&d
\end{pmatrix}\begin{pmatrix}
    1&0\\
    \beta^2&1
\end{pmatrix}=\begin{pmatrix}
    a+\beta^2 b&b\\
    c+\beta^2 d&d
\end{pmatrix}.$$
Then $\beta^2 a=(\beta^2 a+c)-c\in \bbZ$. And similarly, $\beta^2 d\in \bbZ$. In particular, $\mathrm{Tr}(\bar{\Gamma})\subset \frac{1}{\beta^2}\bbZ.$

Let $B\in \bar{\Gamma}.$ If $\mathrm{tr}(B)=\frac{p}{q}$ with $p,q\in \bbZ$, $(p,q)=1$. Then $\mathrm{tr}(B^2)=\mathrm{tr}(B)^2-2$ has a denominator $q^2$. Inductively, $\mathrm{tr}(B^{2^n})$ has a denominator $q^{2^n}$. It follows that $q=1$. In other words, $\mathrm{tr}(B)\in \bbZ.$

$\bar{\Gamma}$ is derived from a quaternion algebra by Theorem~\ref{THM: Fuchsian}. 
\end{proof}
\subsection{Proof of Theorem~\ref{main}}

\begin{proof}
It was shown by Luo and Sarnak in \cite{luo1994number} that if $\Gamma$ is arithmetic, then $\mathrm{Tr}(\Gamma)$ satisfies the B-C property. Furthermore, when $\Gamma$ is derived from a quaternion algebra, $\mathrm{Gap(Tr}(\Gamma))>0$.

On the other hand, if $\mathrm{Tr}(\Gamma)$ satisfies the B-C property, and $\Gamma$ is not uniform, Lemma~\ref{Z} shows that $\bar{\Gamma}$ is derived by quaternion algebra. Hence $\Gamma$ is arithmetic. 

When $\Gamma$ is not derived by quaternion algebra, 
    there must exist an $A\in\Gamma\setminus \bar{\Gamma}$. Left or right multiply by $\begin{pmatrix}
        1&0\\
        \beta^2&1
    \end{pmatrix}$ if necessary, $A$ is in the form  $\begin{pmatrix}
        *&*\\
        c&*
    \end{pmatrix}$ with $c\notin \bbQ$. Hence, $kc+\mathrm{tr}(A)\in \mathrm{Tr}(\Gamma)$ for all $k\in \bbZ$. On the other hand, $2+l\beta^2\in \mathrm{Tr}(\Gamma)$ for all $l\in \bbZ$. 
By Kronecker theorem, $\mathrm{Gap}(\mathrm{Tr}(\Gamma))=0$.
\end{proof}
\section{Klein groups}
The proof of Theorem~\ref{main} for Klein groups follows a similar strategy to the Fuchsian case. 

Let $\Gamma$ be a non-uniform  Kleinian group and $x\in\partial\Gamma$ a cusp point with a cusp subgroup $\Gamma_x$. Taking $g\in \Gamma$ with $gx\neq x$, then $gx$ is a cusp point of $\Gamma$, and $\Gamma_{gx}=g\Gamma_x g^{-1}$. Up to conjugation in $\PSL(2,\mathbb{C})$, we may assume $x=\infty$, $gx=0$, and
$$\setdef{\begin{pmatrix}
    1&k+l\alpha\\
    0&1
\end{pmatrix}}{k,l\in \bfZ}\subset \Gamma_\infty$$
for some $\alpha\notin \mathbb{R}$.

$g\infty=0$ implies $g$ is in the form $\begin{pmatrix}
    0&\frac{1}{\beta}\\
    -\beta&*
\end{pmatrix}$ for some $\beta\in \bbC$. It follows from the fact $\Gamma_0=g\Gamma_\infty g^{-1}$ that
$$\setdef{\begin{pmatrix}
    1&0\\
    \beta^2(k+l\alpha)&1
\end{pmatrix}}{k,l\in \bfZ}\subset \Gamma_0.$$

From now on, we assume $\Gamma$ is a Kleinian lattice containing $\Gamma_0$ and $\Gamma_\infty$ as above. 
\subsection{Prestep: $\alpha$ is a quadratic algebraic number}
Unlike the Fuchsian case, the arithmetic property of $\alpha$ is related to the structure of $\Gamma$. Hence we start the proof with the following lemma.

\begin{lemma}\label{alpha}
If $\mathrm{Tr}(\Gamma)$ has linear growth, $\alpha$ is an algebraic number and $\mathrm{deg}(\alpha)=2$.
\end{lemma}
\begin{proof}
    First, the following identity 
    $$\begin{pmatrix}
        1&m+n\alpha\\
        0&1
    \end{pmatrix}\begin{pmatrix}
        1&0\\
        \beta^2(m'+n'\alpha)&1
    \end{pmatrix}=\begin{pmatrix}
        1+\beta^2(m'+n'\alpha)(m+n\alpha)&m+n\alpha\\
        \beta^2(m'+n'\alpha)&1
    \end{pmatrix}$$
    show that $2+\beta^2(m'+n'\alpha)(m+n\alpha)\in \mathrm{Tr}(\Gamma)$ for all $m$, $n$, $m'$, $n'\in \bbZ$. Since $\mathrm{Tr}(\Gamma)$ has linear growth, the same is true for the set 
    $$\setdef{nn'\alpha^2+(nm'+mn')\alpha+mm'}{m,n,m',n'\in \bbZ}.$$

    For any complex number $x$, consider the two maps $f:\bbZ^4\rightarrow \bbZ^3$ and $g_x:\bbZ^4\rightarrow \mathbb{C}$ given by
    $$f(m,n,m',n')=(nn',nm'+mn',mm'),$$
    $$g_x(m,n,m',n')=nn'x^2+(nm'+mn')x+mm'.$$
Let $R_N$ be the subset of $\bbZ^4$ give by 
$$\setdef{(r_1,r_2,r_3,r_4)}{1\leq r_2\leq r_1\leq N, (r_1,r_2)=1, r_4\leq r_3\leq \frac{N}{r_1}, (r_3,r_4)=1}.$$
 \begin{claim}
     $f$ restricted on $R_N$ is at most 2 to 1, and
     $\lim_{N\rightarrow \infty}\frac{\#R_N}{N^2}=\infty.$
 \end{claim}
 Assume the claim for now.
 Direct computation shows that for all $u\in R_N$, $|g_\alpha (u)|\leq 4N(1+|\alpha|^2)$. Linear growth implies that $\#\setdef{g_\alpha(u)}{u\in R_N}\leq LN^2$ for some $L$. The claim insures two different $u$, $v\in R_N$ for $N$ big enough such that $g_\alpha(u)=g_\alpha(v)$ and $f(u)\neq f(v)$. Then $g_\alpha(u)=g_\alpha(v)$ is a quadratic (possibly linear) polynomial equation for $\alpha$ with integer coefficients. If $\mathrm{deg}(\alpha)=1$, then $\alpha\in\bbQ$, a contradiction. Hence $\mathrm{deg}(\alpha)=2$. 
    
Now we prove the claim. First, we estimate the number of elements in $R_N$. Fix $r_1$, $r_2$ has $\varphi(r_1)$ many choices where $\varphi$ is the Euler's totient function, similar for $r_4$. Hence
$$\#R_N=\sum_{i=1}^N\varphi(i)\sum_{j=1}^{\frac{N}{i}}\varphi(j).$$
Recall in \cite{walfisz1963weylsche}, Walfisz improved Mertens results and get
$$\sum_{1\leq n\leq x}\varphi(n)=\frac{3}{\pi^2}x^2+O(x\log^{0.75}x(\log\log x)^2).$$
Hence
$$\#R_N\geq \kappa\sum_{i=2}^{\frac{N}{2}}\varphi(i)(\frac{N}{i})^2=\kappa N^2\sum_{i=2}^{\frac{N}{2}}\varphi(i)(\frac{1}{i})^2$$
for some $\kappa>0$.

It is well-known (for example, see \cite{rosser1962approximate}) that for $n\geq 2$
$$\varphi(n)>\frac{n}{e^{\gamma}\log\log n+\frac{3}{\log\log n}},$$
where $\gamma$ is the Euler constant. It leads to
$$\#R_N>\kappa N^2\sum_{i=2}^{\frac{N}{2}}\frac{1}{i e^{\gamma}\log\log i+\frac{3i}{\log\log i}}.$$ 
The estimate is a consequence of the divergence of the series. 

\medskip

Second, we need to show $f$ is at most 2 to 1 on $R_N$. Assume $u=(u_1,u_2,u_3,u_4)$ and $v=(v_1,v_2,v_3,v_4)$ are two elements of $R_N$ with $f(u)=f(v)$. Then $g_x(u)=g_x(v)$ for all $x$. In particular, $g_{-\frac{u_1}{u_2}}(u)=0$. Hence 
$v_2(-\frac{u_1}{u_2})+v_1=0$ or $v_3(-\frac{u_1}{u_2})+v_4=0$.

If $v_2(-\frac{u_1}{u_2})+v_1=0$, then $\frac{v_1}{v_2}=\frac{u_1}{u_2}$. By assumption, both numbers are in reduced form, hence $v_1=u_1$, $v_2=u_2$. The first and third coordinates of $f$ tell us that $u=v$.

If $v_3(-\frac{u_1}{u_2})+v_4=0$. Then similar $v_3=u_1$ and $v_4=u_2$. Repeating the same argument for $g_{-\frac{u_3}{u_4}
}$, $v_1=u_3$, $v_2=u_4$.

In summary, $f$ is injective on the diagonal $\{(r_1,r_2,r_1,r_2)\}\cap R_N$ and 2 to 1 otherwise. 
\end{proof}

\subsection{Step one: $\Gamma^{(2)}$ is algebraic.}
A similar method as in the Fuchsian case applies to Kleinian groups.  Recall that $\alpha$ is a degree 2 algebraic number.

Let $x_1,x_2,\cdots,x_n$ be complex numbers. Denote $\bbQ(\alpha)\langle x_1,x_2,\cdots,x_n\rangle$ the $\bbQ(\alpha)$-vector space generated by $x_1$, $x_2,\cdots, x_n$. 
\begin{lemma} If $\mathrm{Tr}(\Gamma)$ has linear growth, then $A^2\in \PSL(2,\bbQ(\alpha))$ for all $A\in \Gamma$.
\end{lemma}
\begin{proof}
Let $A=\begin{pmatrix}
    a&b\\c&d
\end{pmatrix}\in \Gamma$.

First, we assume $c\neq 0$.

From the equation
$$\begin{pmatrix}
  1&0\\
  \beta^2(m+n\alpha)&1
\end{pmatrix}\begin{pmatrix}
    a&b\\
    c&d
\end{pmatrix}\begin{pmatrix}
    1&m'+n'\alpha\\
    0&1
\end{pmatrix}$$
$$=\begin{pmatrix}
    a&a(m'+n'\alpha)+b\\
    \beta^2(m+n\alpha)a+c&\beta^2(m+n\alpha)(m'+n'\alpha)a+(m'+n'\alpha)c+\beta^2(m+n\alpha)b+d
\end{pmatrix},$$
$a+d+\beta^2(m+n\alpha)(m'+n'\alpha)a+(m'+n'\alpha)c+\beta^2(m+n\alpha)b\in \mathrm{Tr}(\Gamma)$. Hence the set
$$\setdef{\beta^2(m+n\alpha)(m'+n'\alpha)a+(m'+n'\alpha)c+\beta^2(m+n\alpha)b}{m,n,m',n'\in\bbZ}$$
satisfies the B-C property.

Denote $\Theta(m, n,m',n')=\beta^2(m+n\alpha)(m'+n'\alpha)a+(m'+n'\alpha)c+\beta^2(m+n\alpha)b$.

\textbf{Step 1.} $\beta^2 a\in \bbQ\langle c,\beta^2b\rangle.$ 
 
Assume $\beta^2 a\notin \bbQ\langle c,\beta^2b\rangle.$ The same argument as in proof of Lemma~\ref{A^2} shows that $\Theta$ is at most 2 to 1 when $(m+n\alpha)(m'+n'\alpha)\neq 0$.

On the other hand, the size of the set 
$$D'_N:=\setdef{(m,n,m',n')}{1\leq |(m+n\alpha)(m'+n'\alpha)|\leq N,m,n,m',n'\in \bbZ}$$
grows as $N^2\ln N$, and $\Theta$ map this set to numbers of complex norm less that $AN(|\beta^2 a|+|\beta^2 b|+|c|)$ where $A=\max\{1, \frac{1}{|m+n\alpha|}|(m,n)\in D_N'\}$. Contrary to linear growth.

\textbf{Step 2.} $\beta^2 b\in \bbQ(\alpha)\langle c\rangle.$ 

There exist $s,t\in \bbQ(\alpha)$ with $\beta^2a=s\beta^2b+tc$. The same proof as in Lemma~\ref{A^2} with $D_N$ replaced by $D'_N$ applies.

\textbf{Step 3.} $A^2\in \PSL(2,\bbQ(\alpha))$.

Step 1 and 2 shows that $\beta^2 a\in \bbQ(\alpha)\langle c\rangle$.
Same argument for  
$$\begin{pmatrix}
  1&m'+n'\alpha\\
  0&1
\end{pmatrix}\begin{pmatrix}
    q&b\\
    c&d
\end{pmatrix}\begin{pmatrix}
    1&0\\
    \beta^2(m+n\alpha)&1
\end{pmatrix}$$
implies $\beta^2 d\in \bbQ(\alpha)\langle c\rangle.$

Consider $\begin{pmatrix}
    1&1\\
    0&1
\end{pmatrix}\begin{pmatrix}
    1&0\\
    \beta^2&1
\end{pmatrix}=\begin{pmatrix}
    1+\beta^2&1\\
    \beta^2&1
\end{pmatrix}\in\Gamma$. Then $\beta^2+\beta^4\in \bbQ(\alpha)\langle \beta^2\rangle$. We conclude that $\beta^2\in \bbQ(\alpha)$. 

Finally, we have $A=cB$ with $B\in \PGL(2,\bbQ(\alpha))$. Taking the determinant, $c^2\in \bbQ(\alpha)$. And $A^2=c^2B^2\in \PSL(2,\bbQ)$.  

If $c=0$, then consider $A\begin{pmatrix}
    1&0\\
    \beta^2&1
\end{pmatrix}=\begin{pmatrix}
    a+\beta^2 b&b\\
    \beta^2 d& d
\end{pmatrix}$. Hence $a,b\in\bbQ(\alpha)\langle d\rangle.$ Now $A^2\in \PSL(2,\bbQ(\alpha))$ follows the same way.
\end{proof}

\subsection{Step two: $\mathrm{Tr}(\Gamma^{(2)})$ is in algebraic integers.}\label{algebraic}
We start with some background on quadratic fields and estimations. 

$\alpha\in \bbQ(\sqrt{-D})$ for some square-free integer $D$. Let $\mathcal{O}$ be the ring of algebraic integers in $\bbQ(\sqrt{-D})$. 
Recall that 
\[
\mathcal{O}=\left\{
\begin{array}{ll}
   \bbZ[\sqrt{D}],  & D \equiv 2,3\ (\textrm{mod}\ 4), \\
    \bbZ[\frac{1+\sqrt{D}}{2}], & D \equiv 1\ (\textrm{mod}\ 4).
\end{array}
\right.
\]
It is clear that there exists $M_1\in \bbN^+$ such that $M_1\mathcal{O}\subset \setdef{m+n\alpha}{m,n\in \bbZ}$. Fix such a $M_1$.

\begin{definition}For any $r,s\in \mathcal{O}$, we denote $(r,s)=1$ if there exist $u'', v''\in \mathcal{O}$ satisfying $u''r+v''s=1.$
\end{definition}

Let $(r,s)=1$. \begin{lemma}
    Assume that there exist $s_1\in \mathcal{O}$ such that the principle ideal $(s_1)$ is primary, and $s\in (s_1)$. If $(r,s)=1$, then there exist $u, v\in \mathcal{O}$ and constant $M_2$ with $us+vr=1$, $|v|\leq M_2|r|$ and $(v,s_1)=1$.
\end{lemma}
\begin{proof}
By definition, there exist $u'', v''\in \mathcal{O}$ satisfying $u''r+v''s=1.$ Taking a representation $v'$ of $v''$ in $\frac{\mathcal{O}}{r\mathcal{O}}$ with least norm. We have  $v''=v'+wr$ with $|v'|\leq \sqrt{1+D^2}|r|$ and $w\in\mathcal {O}$. Then $(u''-ws)r+v's=1$. If $(v',s_1)=1$, take $v=v'$, and $u=u''-ws$. If $(v',s_1)\neq 1$, then there exist $k$ such that $v'^{k}\in (s_1)$. Take $v=v'+r$, $u=u''-ws-s$.

Now,
\[
(v'+r)\sum_{i=1}^k((-1)^{i+1}v'^{2i-1})+(u''-ws-s)s(\sum_{i=1}^k((-1)^{i+1}v'^{2i-2}))=1+(-1)^{k+1}v'^{2k}.
\]
Since $s\in(s_1)$ and $v^{2k}\in (s_1)$, therefore, $(v,s_1)=1$. Taking $M_{2}=\sqrt{1+D^2}+1$, the lemma follows.
\end{proof}
Take $M=\max\{M_1.M_2\}$. Let $\bar{\Gamma}=\Gamma\cap(\PSL(2,\bbQ(\alpha)))$. $\bar{\Gamma}$ is a finite index subgroup of $\Gamma$ since $\Gamma$ is finitely generated and $\Gamma^{(2)}<\bar{\Gamma}$.  Choose $\begin{pmatrix}
    a&b\\
    c&d
\end{pmatrix}\in \bar{\Gamma}.$
$$
\begin{pmatrix}
    1&1\\
    0&1
\end{pmatrix}\begin{pmatrix}
    a&b\\
    c&d
\end{pmatrix}\begin{pmatrix}
    1&-1\\
    0&1
\end{pmatrix}\begin{pmatrix}
    d&-b\\
    -c&a
\end{pmatrix}=\begin{pmatrix}
    1+ac+c^2&1-ac-a^2\\
    c^2&1-ac
\end{pmatrix}\in \bar{\Gamma}.$$
Iteration gives a sequence of elements $A_n\in \bar{\Gamma}$ such that
$A_n=\begin{pmatrix}
   *&*\\
   c^{2^n}&*
\end{pmatrix}$ and $\mathrm{tr}(A_n)=2+c^{2^n}$.
From the identity $\mathrm{tr}(A_n\begin{pmatrix}
    1&m+n\alpha\\
    0&1
\end{pmatrix})=2+(m+1+n\alpha)c^{2^n}$,
 the set
$$\Delta_c:=\setdef{M_1xc^{2^n}}{x\in \mathcal{O}, n\in\bbN}$$
satisfies the B-C property.

\begin{lemma}\label{uv}
    If $c\notin \mathcal{O}$, then $\Delta_c$ has no B-C property. 
\end{lemma}
\begin{proof}
Let $2^sh$ be the order of the ideal class group of $\mathcal{O}$ with $h$ odd. Choose  $t\in \bbN^+$ so that $2^t\equiv 1(\mathrm{mod}\  h)$. We fix this choice once and for all. Since $c\notin \mathcal{O}$, $c^{2^s(2^t-1)}=\frac{p}{q}$ with $(p,q)=1$, $p,q\in \mathcal{O}$, $N_{\bbQ(\alpha)}q\neq 1$. Let $q_1$ be a factor of $q$ such that the principal ideal $(q_1)$ is primary. This is possible since $2^sh$ is the order of the ideal class group.

Follows from $c^{2^{s+t}}=c^{2^s}C^{2^s(2^t-1)}$, there exists a strictly increasing function $g:\bbN^+\rightarrow \bbN^+$ such that $c^{2^{s+kt}}=(\frac{p}{q})^{g(k)}c^{2^s}$.

 Define a function $f$ on $\bbN$ as follows: $f(0)=0$,
$$f(n)=\min\setdef{g(k)}{|q|^{g(k)}\geq 2M^n|c^{2^s}|\prod_{i=o}^{n-1}|p|^{f(i)}}.$$

By Lemma~\ref{uv}, induction starting from $i=n,$  there are $u_i$, $v_i\in\mathcal {O}$  such that
$$u_ip^{f(i)}-v_iv_{i+1}\cdots v_nq^{f_i}=1,$$ $(u_i,q)=1$, $(v_i,q_1)=1$, $|v_i|\leq M|p|^{f(i)}$. 

Let $m_0=M_1\prod_{i=1}^n v_i$, $m_k=M_1u_k\prod_{i=1}^{k-1}v_i$, $1\leq i\leq n$. Then $z_j=m_jc^{2^s}(\frac{p}{q})^{f(i)}\in \Delta_c$ for all $0\leq j\leq n$.

First, since $(u_i,q_1)=1$ and $(v_i,q_1)=1$,  $z_j$, $0\leq i\leq n$  are all different by considering $\bar{z_i}=\frac{z_i}{M_1}$.  

Second, 
\begin{equation*}
    \begin{array}{rl}
         |z_j-z_0|&=|M_1c^{2^s}\prod_{i=1}^{j-1}v_i||u_i(\frac{p}{q})^{f(j)}-\prod_{i=j}^nv_i|  \\
         &=\frac{|M_1c^{2^s}\prod_{i=1}^{j-1}v_i|}{|q|^{f(j)}}|u_ip^{f(i)}+v_iv_{i+1}\cdots v_nq^{f(i)}|\\
         &\leq \frac{M^n|c^{2^s}|\prod_{i=o}^{j-1}|p|^{f(i)}}{|q|^{f(j)}}\\
         &\leq \frac{1}{2}.
    \end{array}
\end{equation*}
It follows that $\Delta_c$ has no $B-C$ property.
\end{proof}
We conclude that $c\in \mathcal{O}$. 

Consider the following two equations:
$$\begin{pmatrix}
    1&0\\
    \beta^2&1
\end{pmatrix}\begin{pmatrix}
    a&b\\
    c&d
\end{pmatrix}=\begin{pmatrix}
    a&b\\
    \beta^2 a+c&\beta^2 b+d
\end{pmatrix},$$
$$\begin{pmatrix}
    a&b\\
    c&d
\end{pmatrix}\begin{pmatrix}
    1&0\\
    \beta^2&1
\end{pmatrix}=\begin{pmatrix}
    a+\beta^2 b&b\\
    c+\beta^2 d&d
\end{pmatrix}.$$
Then $\beta^2 a=(\beta^2 a+c)-c\in \mathcal{O}$, and similarly, $\beta^2 d\in \mathcal{O}$. In particular, $\mathrm{Tr}(\bar{\Gamma})\subset \frac{1}{\beta^2}\mathcal{O}.$ 

Take $S\in \bbN^+$ so that $S\mathrm{Tr}(\bar{\Gamma})\subset \mathcal{O}$. Let $B\in \bar{\Gamma}.$ Since 
$$\mathrm{tr}(B^n)=\mathrm{tr}(B)\mathrm{tr}(B^{n-1})-\mathrm{tr}(B^{n-2}).$$ Induction shows that $S(\mathrm{tr}(B))^n\in \mathcal{O}$.  Taking norm, we have $N_{\bbQ(\alpha)}(\mathrm{tr}(B)))\in \bbZ$. Then by
$$\mathrm{Tr}^n_{\bbQ(\alpha)}(\mathrm{tr}(B))=\mathrm{Tr}_{\bbQ(\alpha)}(\mathrm{tr}(B))\mathrm{Tr}^{n-1}_{\bbQ(\alpha)}(\mathrm{tr}(B))+N_{\bbQ(\alpha)}(\mathrm{tr}(B))\mathrm{Tr}^{n-2}_{\bbQ(\alpha)}(\mathrm{tr}(B)),$$
$S\mathrm{Tr}^n_{\bbQ(\alpha)}(\mathrm{tr}(B))\in\bbZ$. Hence  $\mathrm{Tr}_{\bbQ(\alpha)}(\mathrm{tr}(B))\in\bbZ$.

In other words, $\mathrm{tr}(B)\in \mathcal{O}.$

$\bar{\Gamma}$ is derived from a quaternion algebra by Theorem~\ref{THM: Kleinian}.

\subsection{Proof of Theorem~\ref{main}}
\begin{proof}
  It is not hard to generalize Luo and Sarnak's work in \cite{luo1994number} to Kleinian groups.

On the other hand, if $\mathrm{Tr}(\Gamma)$ satisfies the B-C property, and $\Gamma$ is not uniform, \S~\ref{algebraic} show that $\bar{\Gamma}$ is derived by quaternion algebra. Hence $\Gamma$ is arithmetic. 

When $\Gamma$ is not derived by quaternion algebra, 
    there must exist an $A\in\Gamma\setminus \bar{\Gamma}$. Left or right multiply by $\begin{pmatrix}
        1&0\\
        \beta^2&1
    \end{pmatrix}$ if necessary, $A$ is in the form  $\begin{pmatrix}
        *&*\\
        c&*
    \end{pmatrix}$ with $c\notin \bbQ(\alpha)$. Hence, $kc+\mathrm{tr}(A)\in \mathrm{Tr}(\Gamma)$ for all $k\in M_1\mathcal{O}$. On the other hand, $2+l\beta^2\in \mathrm{Tr}(\Gamma)$ for all $l\in M_1\mathcal{O}$. 
By Kronecker theorem, $\mathrm{Gap}(\mathrm{Tr}(\Gamma))=0$.
\end{proof}

\section{A corollary}
In \cite{lafont2019primitive}, Lafont and McReynolds showed that every noncompact, locally symmetric, arithmetic manifold has arbitrarily long arithmetic progressions
in its primitive length spectrum. This result was extended to every arithmetic locally symmetric orbifold of classical type without Euclidean or compact factors by Miller \cite{miller2016arithmetic}.

We also consider the arithmetic properties of the trace set. Let the \emph{reduced trace set} $\mathrm{Tr_{R}}(\Gamma)$ be the set of traces of non-trivial elements in $\Gamma$.
\begin{corollary}\label{Group}
    Let $\Gamma$ be a lattice of  $\PSL(2,\mathbb{R})$ or $\PSL(2,\mathbb{C})$. If $\mathrm{Tr_R}(\Gamma)$ is closed under subtraction, i.e. $a-b\in \mathrm{Tr_R}(\Gamma)$ for all $a$, $b\in \mathrm{Tr_R}(\Gamma)$, then $\Gamma$ is non-uniform and derived from a quaternion algebra.
\end{corollary}

\begin{proof}
Notice our definition of trace set is invariant under the change of a sign. Since $a$, $b\in \mathrm{Tr}(\Gamma)$ implies $a+b=a-(-b)\in \mathrm{Tr}(\Gamma)$, $\mathrm{Tr}(\Gamma)$ is a group under addition. 

It is clear 
$$\mathrm{Tr}(A^2)=\mathrm{Tr}^2(A)-2.$$
Let $a$, $b\in\mathrm{Tr}(\Gamma)$. Then 
$$(a+b)^2-2+(a-b)^2-2-2(a^2-2)-2(b^2-2)=4\in \mathrm{Tr}(\Gamma),$$
and 
$$4^2-2-3\times 4=2\in\mathrm{Tr}(\Gamma).$$
Hence $\Gamma$ is non-uniform. 

Since $\mathrm{Tr}(\Gamma)$ has no accumulation point, we know that $\mathrm{Gap}(\mathrm{Tr}(\Gamma))>0$. Then $\Gamma$ is derived from a quaternion algebra by Theorem~\ref{main}.
\end{proof}
\begin{remark}
 In \cite{lakeland2017equivalent}, Lakeland has constructed  families of infinitely many lattices with trace  set $\bbZ$ and $2\bbZ$.  
\end{remark}


\begin{bibdiv}
\begin{biblist}
\bib{borel1981commensurability}{article}{
  title={Commensurability classes and volumes of hyperbolic 3-manifolds},
  author={Borel, Armand},
  journal={Annali della Scuola Normale Superiore di Pisa-Classe di Scienze},
  volume={8},
  number={1},
  pages={1--33},
  year={1981}
}

\bib{geninska2008geometric}{article}{
  title={A geometric characterization of arithmetic Fuchsian groups},
  author={Geninska, Slavyana},
  author={Leuzinger, Enrico},
  journal={Duke Mathematical Journal},
  volume={142},
  number={1},
  pages={111--125},
  year={2008},
  publisher={Duke University Press}
}

\bib{lafont2019primitive}{article}{
  title={Primitive geodesic lengths and (almost) arithmetic progressions},
  author={Lafont, J-F}, author={McReynolds, David Ben},
  journal={Publicacions Matem{\`a}tiques},
  volume={63},
  number={1},
  pages={183--218},
  year={2019},
  publisher={Universitat Aut{\`o}noma de Barcelona, Departament de Matem{\`a}tiques}
}

\bib{lakeland2017equivalent}{article}{
  title={Equivalent trace sets for arithmetic Fuchsian groups},
  author={Lakeland, Grant},
  journal={Proceedings of the American Mathematical Society},
  volume={145},
  number={1},
  pages={445--459},
  year={2017}
}

\bib{luo1994number}{article}{title={Number variance for arithmetic hyperbolic surfaces},
  author={Luo, Wenzhi}, author={Sarnak, Peter},
  journal={Communications in mathematical physics},
  volume={161},
  number={2},
  pages={419--432},
  year={1994},
  publisher={Springer}
}

\bib{maclachlan1987commensurability}{inproceedings}{
  title={Commensurability classes of arithmetic Kleinian groups and their Fuchsian subgroups},
  author={Maclachlan, Colin},
  author={Reid, Alan W},
  booktitle={Mathematical Proceedings of the Cambridge Philosophical Society},
  volume={102},
  number={2},
  pages={251--257},
  year={1987},
  organization={Cambridge University Press}
}

\bib{maclachlan2003arithmetic}{book}{
  title={The arithmetic of hyperbolic 3-manifolds},
  author={Maclachlan, Colin},
  author={Reid, Alan W}, 
  volume={219},
  year={2003},
  publisher={Springer}
}

\bib{miller2016arithmetic}{article}{
  title={Arithmetic Progressions in the Primitive Length Spectrum},
  author={Miller, Nicholas},
  journal={arXiv preprint arXiv:1602.01869},
  year={2016}
}

\bib{rosser1962approximate}{article}{
  title={Approximate formulas for some functions of prime numbers},
  author={Rosser, J Barkley},
  author={Schoenfeld, Lowell},
  journal={Illinois Journal of Mathematics},
  volume={6},
  number={1},
  pages={64--94},
  year={1962},
  publisher={Duke University Press}
}

\bib{MR1321639}{incollection}{
    AUTHOR = {Sarnak, Peter},
     TITLE = {Arithmetic quantum chaos},
 BOOKTITLE = {The {S}chur lectures (1992) ({T}el {A}viv)},
    SERIES = {Israel Math. Conf. Proc.},
    VOLUME = {8},
     PAGES = {183--236},
 PUBLISHER = {Bar-Ilan Univ., Ramat Gan},
      YEAR = {1995},
}

\bib{MR1394753}{article}{
    AUTHOR = {Schmutz, Paul},
     TITLE = {Arithmetic groups and the length spectrum of {R}iemann
              surfaces},
   journal={Duke Mathematical Journal},
    VOLUME = {84},
      YEAR = {1996},
    NUMBER = {1},
     PAGES = {199--215},
      ISSN = {0012-7094},
}

\bib{takeuchi1975characterization}{article}{
  title={A characterization of arithmetic Fuchsian groups},
  author={Takeuchi, Kisao},
  journal={Journal of the Mathematical Society of Japan},
  volume={27},
  number={4},
  pages={600--612},
  year={1975},
  publisher={The Mathematical Society of Japan}
}

\bib{walfisz1963weylsche}{article}{
  title={Weylsche Exponentialsummen in der neueren Zahlentheorie},
  author={Walfisz, Arnold},
  journal={VEB Deutscher Verlag der Wissenschaften},
  year={1963}
}

 \end{biblist}
\end{bibdiv}
\end{document}